\documentclass[12pt]{article}
\usepackage{amsmath,amsthm}
\usepackage{amssymb,latexsym}
\usepackage{enumerate}
\usepackage[english]{babel}
\usepackage{graphicx}

\def\Q{{\mathbb Q}}
\def\Z{{\mathbb Z}}

\newtheorem{lemma}{Lemma}
\newtheorem{theorem}[lemma]{Theorem}

\title{
Calculating power integral bases by solving relative Thue equations
}
\author{
Istv\'{a}n Ga\'{a}l\thanks{
        This work was partially supported by the European Union 
        and the European Social Fund through project Supercomputer, 
        the national virtual lab (grant no.: TAMOP-4.2.2.C-11/1/KONV-2012-0010)
        and also supported in part by K67580 and K75566 from the
        Hungarian National Foundation for Scientific Research,
         },\; 
L\'aszl\'o Remete and T\'\i mea Szab\'o
\\ \\
University of Debrecen, Mathematical Institute \\
H--4010 Debrecen Pf.12., Hungary \\
e--mail: igaal@science.unideb.hu, remetel42@gmail.com, \\ szabo.timea@science.unideb.hu
}

\begin{document}

\maketitle
\thispagestyle{empty}

\renewcommand{\thefootnote}{}



\renewcommand{\thefootnote}{\arabic{footnote}}
\setcounter{footnote}{0}

\begin{abstract}
In our recent paper \cite{girelthue} we gave an efficient algorithm to calculate "small"
solutions of relative Thue equations (where "small" means an upper bound of type $10^{500}$
for the sizes of solutions). Here we apply this algorithm to calculating power integral bases
in sextic fields with an imaginary quadratic subfield and to calculating relative
power integral bases in pure quartic extensions of imaginary quadratic fields. 
In both cases the crucial point of the calculation is the resolution of a relative
Thue equation. We produce numerical data that were not known before.
 \end{abstract}

\section{Introduction}

Calculating power integral bases is a classical field of algebraic number theory
(c.f. \cite{book}). Considering several types of number fields we have seen that
this problem often leads to the resoltion of various types of Thue equations
\cite{thueapple}, \cite{compos}, \cite{degnine}, \cite{relcubic},
\cite{gpezaz}, \cite{relquartic}.

We often used the method of A.Peth\H o \cite{pet}, based on the continued fraction
algorithm, which gave an efficient way to calculate "small" solutions of Thue equations.
"Small" yields here an upper bound, say $10^{500}$, for the absolutie values
of the solutions. This was very much faster than the complete resolution of the equation, 
and gave all solutions with very high probability,
certainly all that can be used in practice. 
(Our experience shows that such equations usually only have a few and rather small
solutions.) It made possible to get an overview
on the solutions of a huge number of equations. Applying it to problems on
power integral bases, we got an overview about the existence of
power integral bases of a huge number of
fields (cf. e.g. \cite{gaalremete}).

We have also learnt that calculating power integral bases in higher degree number
fields having certain subfields often leads to relative Thue equations
\cite{relcubic}, \cite{gpezaz}, \cite{relquartic}. Although there is an algorithm for the
complete resolution of Thue equations \cite{relthue}, an analogue of Peth\H o's
fast algorihtm \cite{pet} in the relative case was missing. Recently the first author \cite{girelthue}
developed such a fast algorithm to calculate "small" solutions (e.g. with sizes less than
$10^{500}$) of relative Thue equations. The algorithm is based on the LLL reduction algorithm
\cite{lll} as one could expect. Since in higher degree number fields even the
calculation of basic field data (integral basis, fundamental units) can become
a hard and time consuming problem, this algorithm seems to have several useful
applications.

In this paper we present two applications. First we calculate power integral bases
in sextic fields with an imaginary subfield. This type of problem we studied already
in \cite{gpezaz}, our purpose to reconsider it is to compare CPU times and to extend
the short list of results we got in \cite{gpezaz}.

Second, we consider pure quartic extensions of imaginary quadratic fields and
calculate relative power integral bases in these relative quartic extensions. This 
calculation generalizes our recent results \cite{gaalremete} on power integral
bases in pure quartic fields.

\section{Basic concepts of power integral bases}

Let $K$ be an algebraic number field of degree $n$ 
with ring of integers $\Z_K$. 
An important task in algebraic number theory is to
decide if $K$ has {\em power integral bases}, 
that is integral bases of the form $\{ 1,\gamma,\gamma^2,\ldots,\gamma^{n-1}\}$,
and to determine all generators $\gamma$ of power integral bases  (cf. \cite{book}).
In general, if $\{1,\omega_2,\ldots,\omega_n\}$ is an integral basis of $K$, then
the discriminant of the linear form $\omega_2 X_2+\ldots+\omega_n X_n$
can be written as
\[
D_{K/Q}(\omega_2 X_2+\ldots+\omega_n X_n)=(I(X_2,\ldots,X_n))^2 D_{K}
\]
where $D_K$ denotes the discriminant of the field $K$,
and $I(X_2,\ldots,X_n)$ is the {\em index form} 
corresponding to the above integral basis.
As is known, for any primitive element $\gamma=x_1+\omega_2 x_2 +\ldots +\omega_n x_n\in\Z_K$
(that is $K=\Q(\gamma)$) we have
\[
I(\gamma)=(\Z_K^+:\Z[\gamma]^+)=|I(x_2,\ldots,x_n)|,
\]
where the index of the additive groups of the corresponding rings are taken.
$I(\gamma)$ is the {\em index of} $\gamma$, which does not depend on $x_1$.
Therefore the element $\gamma=x_1+\omega_2 x_2 +\ldots +\omega_n x_n$
generates a power integral basis of $K$ if and only if
$x_1\in\Z$ and $(x_2,\ldots,x_n)\in \Z^{n-1}$ is a solution of the
{\em index form equation} 
\begin{equation}
I(x_2,\ldots,x_n)=\pm 1 \;\;\; {\rm with}\;\;\; x_2,\ldots,x_n\in\Z .
\label{ifeq}
\end{equation}
This is the way the problem of power integral bases reduces to the
resolution of diophantine equations. 

There is an extensive literature of index form equations and power integral bases
(for a summary cf. \cite{book}).

\section{Sextic fields with an imaginary quadratic subfield}

In this section we recall the method of \cite{gpezaz} based on the input data
of M.Olivier \cite{olivier}, calculating power integral bases in sextic fields
with an imaginary quadratic subfield. In \cite{gpezaz} we calculated all 
generators of power integral bases in 25 sextic fields of this type. 
Here we only give the most important steps of the calculation, our point is to 
compare CPU times with that of \cite{gpezaz} and to extend 
considerably the list of \cite{gpezaz} by applying the method of \cite{girelthue}.

In the table of \cite{olivier} the discriminant $D_K$ of the sextic field $K$,
the discriminant $D_L$ of its imaginary quadratic subfield and the cubic minimal 
polynomial $f(x)$ of $\vartheta$ generating $K$ over $L$ is given.

Let $\{ 1,\omega\}$ be the integral basis of $L$. The element $\vartheta$ is choosen
(see \cite{gpezaz}) to have relative index 1 over $L$ so that all
integral elements of $K$ can be written in the form
\begin{equation}
\alpha=x_0+x_1\vartheta+x_2\vartheta^2+y_0\omega+y_1\omega\vartheta+y_2\omega\vartheta^2
\label{alfa}
\end{equation}
with $x_i,y_i\in\Z\;\; (i=0,1,2)$. Let $\varrho=-\vartheta^{(1)}-\vartheta^{(2)}$ with two
distinct roots of $f(x)$. The element $\alpha$ of (\ref{alfa}) generates a power integral basis
of $K$ if and only if the quadratic integers $X=x_1+\omega y_1,\; Y=x_2+\omega y_2$ satisfy
the relative Thue equation
\begin{equation}
N_{K/L}(X-\varrho Y)= \nu \;\; {\text in}\; X,Y\in\Z_L
\label{relthue}
\end{equation}
(with a unit $\nu$ of $L$) and $x_1,x_2,y_0,y_1,y_2$ is a solution of a 
degree 9 polynomial equation
\begin{equation}
F(x_1,x_2,y_0,y_1,y_2)=\pm 1 \;\; {\text in}\; x_1,x_2,y_0,y_1,y_2\in\Z
\label{degkilenc}
\end{equation}
(for the construction of $F$ see \cite{gpezaz}). Therefore in order to determine all
generators of power integral bases (\ref{alfa}) of $K$ one has to determine
$x_1,x_2,y_1,y_2$ by solving the relative Thue equation (\ref{relthue}) 
and to calculate $y_0$ from (\ref{degkilenc})
($x_0\in\Z$ is arbitrary).

In \cite{gpezaz} for the complete resolution of the relative Thue equation one had to 
calculate the basic data (integral basis, fundamental units) of $K$. 
We represented $X-\varrho Y$ as a power product of the fundamental units of $K$ and
the calculation was focused on the unknown exponents. We used Baker's
method, reduction of the bounds and then enumerated the small exponents.
The CPU time was about 20 minutes per example.

In our present calculation we compute the solutions of (\ref{relthue}) with
\begin{equation}
\max(\overline{|X|},\overline{|Y|})<10^{250}
\label{cccc}
\end{equation}
($\overline{|\gamma|}$ denotes the {\it size} 
of $\gamma$, that is the maximum absolute value of its conjugates). 
This  yields that we calculate all generators of power integral bases of $K$
with $\max(|x_1|,|x_2|,|y_1|,|y_2|)<C$ where $C$ is of magnitude $10^{250}$
(the exact value of $C$ can be easily calculeted in view of $\omega$).
Our list of solutions contains all solutions with very high probability, 
certainly all that can be used in practical calculations. 
Our calculation is focused on the solutions themselves, we do not need any 
additional calculation of the basic field data.
The running time of this method was about 2-5 minutes per example, therefore 
it was at least 5-10 time faster
and made possible to list the results of the first 100 number fields of
smallest discriminant in absolute value. 

In the following table for brevity we only list the coordinates $(x_1,x_2,y_1,y_2,y_0)$
of those generators of power integral bases of $K$
(satisfying (\ref{cccc})) which have one or more coordinates $\geq 3$ in absolute value.
We indicate if there are no solutions (no generators of power integral bases) at all.
In all lines where only the field data appear we mean that
all solutions have coordinates $\leq 2$ in absolute value.
In the lines where solutions appear we skiped 
the solutions with coordinates $\leq 2$ in absolute value.
The complete list can be fould at 
http://www.math.unideb.hu/~igaal/ .
\\ \\
{\scriptsize
$D_K=-9747$, $\omega=(1+i\sqrt{3})/2$, $f(x)=x^3-(1+\omega)x^2+\omega x+(1-\omega)$ \\
$D_K=-10816$, $\omega=i$, $f(x)=x^3-(1+\omega)x^2+5\omega x-(1+4\omega)$,
(0 2 -3 0 6), (0 2 -3 0 7), (0 1 -2 0 3), \\
(0 1 -1 0 3), (1 -1 2 1 -4),(1 -1 2 1 -3), (1 -2 3 1 -7), 
(1 -2 3 1 -6), (2 2 -3 1 6), (2 1 -1 1 2),(3 1 -1 2 2), \\
(3 0 0 2 -1), (3 -1 2 2 -5), (3 -1 2 2 -4), (5 2 -3 3 4), (5 2 -2 3 5)\\
$D_K=-11691$, $\omega=(1+i\sqrt{3})/2$, $f(x)=x^3-(1+\omega)x^2+(-2+2\omega)x+1$,
(2 -2 0 1 -3)\\ 
$D_K=-12167$, $\omega=(1+i\sqrt{23})/2$, $f(x)=x^3-(1+\omega)x^2+(-2+\omega)x+1$ \\
$D_K=-14283$, $\omega=(1+i\sqrt{3})/2$, $f(x)=x^3+(1-\omega)x-1$,
(2 -3 -3 0 4), (3 3 -2 -3 -4)\\
$D_K=-16551$, $\omega=(1+i\sqrt{3})/2$, $f(x)=x^3-(1+\omega)x^2+2x+(-1+\omega)$,
(1 1 -4 1 3), (3 0 -3 2 1)\\ 
$D_K=-16807$, $\omega=(1+i\sqrt{7})/2$, $f(x)=x^3-\omega x^2+(-1+\omega)x+1$ \\
$D_K=-19683$, $\omega=(1+i\sqrt{3})/2$, $f(x)=x^3+(-1+\omega)$ \\
$D_K=-21168$, $\omega=(1+i\sqrt{3})/2$, $f(x)=x^3-x^2+(1-2\omega)x+1$ \\
$D_K=-21296$, $\omega=(1+i\sqrt{11})/2$, $f(x)=x^3-\omega x^2+(-1+\omega)x+1$ \\
$D_K=-22592$, $\omega=i$, $f(x)=x^3-(1+\omega)x^2+(1+2\omega)x-\omega$,
(1 2 -3 3 6)\\ 
$D_K=-22707$, $\omega=(1+i\sqrt{3})/2$, $f(x)=x^3-(1+\omega)x^2+2\omega x+(1-2\omega)$ \\
$D_K=-23031$, $\omega=(1+i\sqrt{3})/2$, $f(x)=x^3-x^2+(-1+\omega)$ \\
$D_K=-24003$, $\omega=(1+i\sqrt{3})/2$, $f(x)=x^3-x^2-x+(1-\omega)$ \\
$D_K=-25947$, $\omega=(1+i\sqrt{3})/2$, $f(x)=x^3+x+1$ \\
$D_K=-29791$, $\omega=(1+i\sqrt{31})/2$, $f(x)=x^3-(1+\omega)x^2+(-2+\omega)x+1$ \\
$D_K=-30976$, $\omega=i$, $f(x)=x^3-x^2+(2-\omega)x-1$ \\
$D_K=-31347$, $\omega=(1+i\sqrt{3})/2$, $f(x)=x^3-(1+\omega)x^2+3\omega x-\omega$ \\
$D_K=-33856$, $\omega=i$, $f(x)=x^3+x-\omega$ \\
$D_K=-34371$, $\omega=(1+i\sqrt{3})/2$, $f(x)=x^3-(1+\omega)x^2+(-1+4\omega)x+(2-\omega)$,
(1 -1 0 0 -3), (2 -1 -3 2 2),\\ (3 -2 -4 2 0)\\ 
$D_K=-34992$, $\omega=(1+i\sqrt{3})/2$, $f(x)=x^3-(1+\omega)x^2+3\omega x+(1-2\omega)$ \\
$D_K=-36963$, $\omega=(1+i\sqrt{3})/2$, $f(x)=x^3-(1+\omega)x^2+\omega x+(-1+\omega)$ \\
$D_K=-40203$, $\omega=(1+i\sqrt{3})/2$, $f(x)=x^3-(1+\omega)x^2+(-2+3\omega)x+(2-\omega)$ \\
$D_K=-41472$, $\omega=i\sqrt{2}$, $f(x)=x^3+(1-\omega)x-1$ \\
$D_K=-41823$, $\omega=(1+i\sqrt{3})/2$, $f(x)=x^3-x^2+(5-5\omega)x+(-6+2\omega)$,
(0 -2 -2 1 7), (0 -1 -1 1 3),\\ (1 1 0 -1 -3), (2 -1 -2 -1 5), (3 1 0 -3 -1)\\ 
$D_K=-44496$, $\omega=(1+i\sqrt{3})/2$, $f(x)=x^3-(1+\omega)x^2+2x-2$, no solutions\\ 
$D_K=-47680$, $\omega=i$, $f(x)=x^3-\omega x-1$, (3 3 1 0 -1) \\
$D_K=-47979$, $\omega=(1+i\sqrt{3})/2$, $f(x)=x^3-x^2-2\omega x+(-1+2\omega)$ \\
$D_K=-49408$, $\omega=i$, $f(x)=x^3-x^2+x+\omega$ \\
$D_K=-50139$, $\omega=(1+i\sqrt{3})/2$, $f(x)=x^3-x^2-x+(-1+\omega)$,
(4 -2 -3 2 -1) \\
$D_K=-52272$, $\omega=(1+i\sqrt{3})/2$, $f(x)=x^3-(1+\omega)x^2+(4-\omega)x+(-3-2\omega)$,
(0 -1 -1 1 3), (1 -2 -1 1 5)\\
$D_K=-53568$, $\omega=(1+i\sqrt{3})/2$, $f(x)=x^3-x^2-\omega x-\omega$,
no solutions\\
$D_K=-53824$, $\omega=i$, $f(x)=x^3-(1+\omega)x^2+(2+2\omega)x-1$,
no solutions\\
$D_K=-54675$, $\omega=(1+i\sqrt{3})/2$, $f(x)=x^3-(1+\omega)x^2+(2-\omega)x+\omega$ \\
$D_K=-57591$, $\omega=(1+i\sqrt{3})/2$, $f(x)=x^3-(1+\omega)x^2+(3-3\omega)x+1$,
 (1 -2 4 0 3), (2 0 -3 1 1),\\(3 -1 -1 1 1) \\
$D_K=-59648$, $\omega=i$, $f(x)=x^3-(1+\omega)x^2+(-1+\omega)x+(1+\omega)$ \\
$D_K=-59967$, $\omega=(1+i\sqrt{3})/2$, $f(x)=x^3-(1+2\omega)x-\omega$ \\
$D_K=-60992$, $\omega=i$, $f(x)=x^3-x^2-(1+\omega)x+1$,
(8 -4 1 0 1) \\
$D_K=-61504$, $\omega=i$, $f(x)=x^3+x-1$ \\
$D_K=-64827$, $\omega=(1+i\sqrt{3})/2$, $f(x)=x^3-x^2-2x+1$,
no solutions \\
$D_K=-65600$, $\omega=i$, $f(x)=x^3-(1+\omega)x^2-\omega x-1$ \\
$D_K=-70659$, $\omega=(1+i\sqrt{3})/2$, $f(x)=x^3-2x+(1-\omega)$ \\
$D_K=-72716$, $\omega=(1+i\sqrt{7})/2$, $f(x)=x^3-(1+\omega)x^2+(-1+2\omega)x+1$,
(3 -2 0 0 -1) \\
$D_K=-73008$, $\omega=(1+i\sqrt{3})/2$, $f(x)=x^3-x^2+(3-2\omega)x-1$ \\
$D_K=-73467$, $\omega=(1+i\sqrt{3})/2$, $f(x)=x^3-x^2-2\omega x+1$ \\
$D_K=-82496$, $\omega=i$, $f(x)=x^3-(1+\omega)x^2+1$ \\
$D_K=-82971$, $\omega=(1+i\sqrt{3})/2$, $f(x)=x^3-(1+\omega)x^2+(2+\omega)x+(-2+\omega)$ \\
$D_K=-85131$, $\omega=(1+i\sqrt{3})/2$, $f(x)=x^3-x^2+(3-2\omega)x+(-1+\omega)$,
 (1 -1 -1 1 3)\\
$D_K=-86528$, $\omega=i\sqrt{2}$, $f(x)=x^3-\omega x^2-\omega x-1$ \\
$D_K=-87616$, $\omega=i$, $f(x)=x^3-(1+\omega)x^2+(-2+\omega)x+1$ \\
$D_K=-87831$, $\omega=(1+i\sqrt{3})/2$, $f(x)=x^3-(1+\omega)x^2+2\omega x-(1+\omega)$ \\
$D_K=-91719$, $\omega=(1+i\sqrt{3})/2$, $f(x)=x^3-(1+\omega)x^2+(-1+4\omega)x-2\omega$,
 (2 -2 1 1 -5) \\
$D_K=-92416$, $\omega=i$, $f(x)=x^3-(1+\omega)x^2+(-1+\omega)$,
no solutions \\
$D_K=-93987$, $\omega=(1+i\sqrt{3})/2$, $f(x)=x^3-2\omega x+1$,
no solutions\\
$D_K=-94311$, $\omega=(1+i\sqrt{3})/2$, $f(x)=x^3-x^2-3\omega x+(-1+4\omega)$,
(1 0 -1 -1 3) \\
$D_K=-95607$, $\omega=(1+i\sqrt{3})/2$, $f(x)=x^3-x^2+2x+\omega$ \\
$D_K=-96512$, $\omega=i$, $f(x)=x^3-x^2-x-\omega$ \\
$D_K=-96579$, $\omega=(1+i\sqrt{3})/2$, $f(x)=x^3-x^2+(1-\omega)x+(-2+\omega)$ \\
$D_K=-96832$, $\omega=i$, $f(x)=x^3-(1+\omega)x^2+\omega x+\omega$ \\
$D_K=-103383$, $\omega=(1+i\sqrt{3})/2$, $f(x)=x^3-x^2+(3-3\omega)x+(-3+2\omega)$,
 (1 -1 -1 0 3) \\
$D_K=-104112$, $\omega=(1+i\sqrt{3})/2$, $f(x)=x^3-(1+\omega)x^2+(-2+3\omega)x+(1-2\omega)$ \\
$D_K=-104571$, $\omega=(1+i\sqrt{3})/2$, $f(x)=x^3-x^2-(1+2\omega)x+3\omega$ \\
$D_K=-106560$, $\omega=i$, $f(x)=x^3-(1+\omega)x^2-x-1$ \\
$D_K=-107163$, $\omega=(1+i\sqrt{3})/2$, $f(x)=x^3-(1+\omega)x^2+(1-2\omega)x+(-2+3\omega)$,
 (0 1 0 0 -3) \\
$D_K=-107811$, $\omega=(1+i\sqrt{11})/2$, $f(x)=x^3-\omega x^2+(-3+\omega)x+\omega$ \\
$D_K=-108459$, $\omega=(1+i\sqrt{3})/2$, $f(x)=x^3-x^2+(3-3\omega)x+(-2+\omega)$ \\
$D_K=-108544$, $\omega=i$, $f(x)=x^3-(1+\omega)x^2+(-2+\omega)x+(1+\omega)$ \\
$D_K=-108731$, $\omega=(1+i\sqrt{7})/2$, $f(x)=x^3-(1+\omega)x^2-x+\omega$ \\
$D_K=-108800$, $\omega=i$, $f(x)=x^3-x^2+(1-2\omega)x+\omega$ \\
$D_K=-109539$, $\omega=(1+i\sqrt{3})/2$, $f(x)=x^3-x^2+3x-\omega$,
 (1 -1 -1 1 3)\\
$D_K=-109744$, $\omega=(1+i\sqrt{19})/2$, $f(x)=x^3-(1+\omega)x^2+(-2+\omega)x+1$ \\
$D_K=-110079$, $\omega=(1+i\sqrt{3})/2$, $f(x)=x^3+(1-2\omega)x+(2-\omega)$ \\
$D_K=-110144$, $\omega=i$, $f(x)=x^3-(1+\omega)x^2+(-3-\omega)x+(2+3\omega)$,
(0 2 -1 0 -3), (1 1 0 1 -3) \\
$D_K=-112192$, $\omega=i$, $f(x)=x^3-x^2+(-2-3\omega)x-2\omega$,
(1 -1 2 -1 3),(4 -2 0 1 1) \\
$D_K=-114399$, $\omega=(1+i\sqrt{3})/2$, $f(x)=x^3+(4-\omega)x+(1-3\omega)$,
(2 -3 -3 0 4) \\
$D_K=-116800$, $\omega=i$, $f(x)=x^3+(1-3\omega)x+(-2+\omega)$ \\
$D_K=-117207$, $\omega=(1+i\sqrt{3})/2$, $f(x)=x^3-x^2-2x+(1+\omega)$ \\
$D_K=-118287$, $\omega=(1+i\sqrt{3})/2$, $f(x)=x^3+(-2-\omega)x+\omega$,
(1 1 -1 -3 4)\\
$D_K=-122256$, $\omega=(1+i\sqrt{3})/2$, $f(x)=x^3-(1+\omega)x^2-x+(1-\omega)$,
no solutions \\
$D_K=-124848$, $\omega=(1+i\sqrt{3})/2$, $f(x)=x^3-4\omega x+(-2+4\omega)$,
(0 -1 -1 0 3), (0 0 -1 -1 3), (1 0 -4 -3 7), (1 3 3 0 -7)\\
$D_K=-129088$, $\omega=i$, $f(x)=x^3-x^2+(1-2\omega)x+(1+\omega)$ \\
$D_K=-130032$, $\omega=(1+i\sqrt{3})/2$, $f(x)=x^3-(1+\omega)x^2+(-5+5\omega)x+(6-5\omega)$,
no solutions \\
$D_K=-130304$, $\omega=i$, $f(x)=x^3-(1+\omega)x^2+(-3+2\omega)x+2\omega$,
(4 -2 1 0 -3)\\
$D_K=-131787$, $\omega=(1+i\sqrt{3})/2$, $f(x)=x^3-x^2+(2-3\omega)x+(-1+\omega)$ \\
$D_K=-133407$, $\omega=(1+i\sqrt{3})/2$, $f(x)=x^3-(1+\omega)x^2-(1+\omega)x-1$,
(4 -2 -3 2 0) \\
$D_K=-133839$, $\omega=(1+i\sqrt{3})/2$, $f(x)=x^3-x^2-(1+\omega)x+2$,
 (1 0 -1 2 -5)\\
$D_K=-134811$, $\omega=(1+i\sqrt{3})/2$, $f(x)=x^3+(2-2\omega)x-(1+\omega)$ \\
$D_K=-137200$, $\omega=(1+i\sqrt{7})/2$, $f(x)=x^3-\omega x^2+(-2+\omega)x+\omega$ \\
$D_K=-137403$, $\omega=(1+i\sqrt{3})/2$, $f(x)=x^3-(1+\omega)x^2+(-2+3\omega)x+(2-3\omega)$,
no solutions \\
$D_K=-139023$, $\omega=(1+i\sqrt{3})/2$, $f(x)=x^3-(1+\omega)x^2+(1+2\omega)x-2$ \\
$D_K=-139520$, $\omega=i$, $f(x)=x^3+(1-\omega)x-(1+\omega)$ \\
$D_K=-139968$, $\omega=(1+i\sqrt{3})/2$, $f(x)=x^3-(1+\omega)x^2-(1+\omega)x+\omega$,
no solutions \\
$D_K=-141939$, $\omega=(1+i\sqrt{3})/2$, $f(x)=x^3-(1+\omega)x^2-x-\omega$ \\
$D_K=-143872$, $\omega=i\sqrt{2}$, $f(x)=x^3-\omega x^2-(1+\omega)x-1$ \\
$D_K=-143883$, $\omega=(1+i\sqrt{3})/2$, $f(x)=x^3-(1+\omega)x^2+\omega x-(1+\omega)$,
no solutions \\
$D_K=-144207$, $\omega=(1+i\sqrt{3})/2$, $f(x)=x^3-(1+\omega)x^2+(2-2\omega)x+(-2+\omega)$,
(2 0 -3 1 1) \\
$D_K=-144448$, $\omega=i$, $f(x)=x^3-(1+\omega)x^2-(1+\omega)x-1$ \\
$D_K=-147008$, $\omega=i$, $f(x)=x^3-(1+\omega)x^2-(1+2\omega)x+1$ \\
$D_K=-147520$, $\omega=i$, $f(x)=x^3-(1+\omega)x^2+(-2+\omega)x+(2+\omega)$ \\
$D_K=-149283$, $\omega=(1+i\sqrt{3})/2$, $f(x)=x^3-x^2+(3-\omega)x+(-2+\omega)$ \\
}
\\ \\

\section{Relative power integral bases of pure quartic extensions of imaginary quadratic fields}

In a recent paper \cite{gaalremete} we considered power interal bases in pure quartic fields.
Using a general result of I.Ga\'al, A.Peth\H o and M.Pohst \cite{gppsim} on quartic fields
this problem could be reduced to the resolution of Thue equations. For pure quartic fields
this Thue equation happend to be a binomial Thue equation of type $x^4-m y^4=\pm 1$.
Using the fast algorithm of A.Peth\H o \cite{pet} we calculated the solutions with
$\max(|x|,|y|)<10^{500}$ of these binomial Thue equations up to $m<10^7$ and used these
data for listing all generators of power integral bases of the pure quartic field
$K=\Q({}^4\sqrt{m})$ with coefficients in absolute values less than $10^{1000}$.

As we shall see in the following, the generalization \cite{relquartic} of \cite{gppsim}
can be applied to calculate relative power integral bases of quartic extensions.
Moreover, if we consider pure quartic extensions $K=L({}^4\sqrt{m})$ of imaginary
quadratic fields $L$, then we obtain a binomial relative Thue equation for the
resolution of which we can efficiently use our algorithm \cite{girelthue}.

To fix our notation we let $d>1$ be a square free integer and $L=\Q(i\sqrt{d})$. 
We assume $-d\equiv 1 (\bmod \ 4)$ therefore $D_L=-d$ and an integer basis of $L$
is $\{1,\omega\}$ with $\omega=(1+i\sqrt{d})/2$. 

Let $m$ be a positive integer and consider the pure quartic field $M=\Q({}^4\sqrt{m})$.
Set $\xi={}^4\sqrt{m}$. Assume that $m$ is square free, $m\neq\pm 1$ and
$m\equiv 2,3 \;(\bmod \ 4)$. According to 
A.Hameed, T.Nakahara, S.M.Husnine and S.Ahmad \cite{hnha}
$\{1,\xi, \xi^2, \xi^3\}$ is an integer basis in $M$ 
with discriminant $D_M=-256m^3$.

As it is well known (\cite{nark}, \cite{compos}), if $(D_L,D_M)=1$ then
\[
\{1,\xi, \xi^2, \xi^3,\omega,\omega\xi, \omega\xi^2, \omega\xi^3\}
\]
is an integral basis in the composite field $K=LM=\Q(i\sqrt{d},{}^4\sqrt{m})$.
To use this property we shall assume that $-d\equiv 1 (\bmod \ 4)$, 
$1<m\equiv 2,3 \;(\bmod \ 4)$ and $(d,m)=1$.
Moreover, since our arguments use that $\Z_L$ has unique factorization,
we assume that $-d$ is any of -3,-7,-11,-19,-43,-67,-163 (cf. \cite{baker}).

Then we can write any $\alpha\in\Z_K$ in the form
\begin{equation}
\alpha=H+X\xi+Y\xi^2+Z\xi^3
\label{aaa}
\end{equation}
with $H,X,Y,Z\in\Z_L$. We apply now the main result of 
I.Ga\'al, A.Peth\H o and M.Pohst \cite{relquartic} to our relative
quartic extension $K/L$. 

\begin{lemma} 
Let
\[
F(U,V)=U(U^2-4mV^2),\;\; Q_1=X^2-mZ^2,\;\; Q_2=Y^2-XZ.
\]
If $\alpha$ of (\ref{aaa}) generates a relative power integral 
basis in $K$ over $L$ then there are $U,V\in\Z_L$ such that 
\[
N_{K/L}(F(U,V))=\pm 1
\]
with
\[
Q_1(X,Y,Z)=U,\;\; Q_2(X,Y,Z)=V.
\]
\label{lemma1}
\end{lemma}

\noindent
{\bf Proof.} This is a direct consequence of \cite{relquartic} using that $\xi$
has relative defining polynomial $x^4-m$ over $L$. $\Box $

\subsection{From index equations to binomial Thue equations}

In this section we show that using Lemma \ref{lemma1}, determining elements of $\Z_K$
generating relative power integral bases over $L$ can be reduced to solving
quartic relative binomial Thue equations over $L$.

\begin{theorem}
If $\alpha=H+X\xi+Y\xi^2+Z\xi^3$ of (\ref{aaa}) generates a relative power
integral basis of $K$ over $L$, then there exist a solution $(X_0,Y_0)\in\Z_L^2$
of the relative binomial Thue equation
\[
X_0^4-mY_0^4=\zeta
\]
(with a unit $\zeta\in L$) such that 
\begin{equation}
X=X_0^2\varepsilon_0,\;\;
Y=\pm X_0Y_0\varepsilon_0,\;\;
Z=Y_0^2\varepsilon_0.
\label{xyz}
\end{equation}
\label{thbb}
\end{theorem}

\noindent
{\bf Proof}.
Assume that $\alpha=H+X\xi+Y\xi^2+Z\xi^3$ of (\ref{aaa}) generates a relative power integral 
basis of $K$ over $L$. 
The first equation of Lemma (\ref{lemma1}) implies that 
\[
U(U^2-4mV^2)=\nu
\]
with a unit $\nu$ in $L$. Moreover $U=\eta$ is also a unit in $L$. Hence
\[
V^2=\frac{\eta^2-\nu/\eta}{4m}
\]
which implies that the only possible value of $V$ is $V=0$.

We utilize now the equations of Lemma \ref{lemma1} concerning $Q_1$ and $Q_2$.
Equation $Q_1(X,Y,Z)=X^2-mZ^2=U=\eta$ implies (by unique factorization in $L$,
being valid by the restricted values of $d$)
that $X,Z$ are coprime in $\Z_L$.
From equation $Q_2(X,Y,Z)=Y^2-XZ=V=0$ we confer $Y^2=XZ$ whence
(using that $X,Z$ are coprime) $X=A^2\varepsilon_x, Z=B^2\varepsilon_y$
with certain $A,B\in\Z_L$ and units $\varepsilon_x,\varepsilon_y\in\Z_L$
such that their product is a square $\varepsilon_x\varepsilon_y=\varepsilon^2$.
Substituting these into equation $Q_1(X,Y,Z)=X^2-mZ^2=U=\eta$ we obtain
\[
A^4\varepsilon_x^2-m B^4\varepsilon_y^2=\eta
\]
whence
\[
(A\varepsilon_x)^4-m B^4\varepsilon_x^2\varepsilon_y^2=\eta\varepsilon_x^2
\]
finally by $\varepsilon_x\varepsilon_y=\varepsilon^2$ we get
\[
(A\varepsilon_x)^4-m (B\varepsilon)^4=\eta\varepsilon_x^2.
\]
This way we arrive at the relative binomial Thue equation
\begin{equation}
X_0^4-mY_0^4=\zeta \;\; {\rm in}\;\; X_0,Y_0\in\Z_L
\label{rrr}
\end{equation}
where $\zeta$ is a unit in $\Z_L$.

Collecting the relations between $X,Y,Z$ and $X_0,Y_0$ we obtain
\[
X=A^2\varepsilon_x=X_0^2\varepsilon_x^{-1},
\]
\[
Z=B^2\varepsilon_y=\frac{Y_0^2}{\varepsilon^2}\varepsilon_y=Y_0^2\varepsilon_x^{-1},
\]
finally
\[
Y=\pm \sqrt{XZ}=\pm X_0Y_0\varepsilon_x^{-1}.
\]
For any $\varepsilon_x$ we can obviously choose an $\varepsilon_y$ such that
their product is a square. As we see, finally $X,Y,Z$ are only related with
$X_0,Y_0$ by means of $\varepsilon_x$, therefore we obtain the statement of Theorem \ref{thbb}
by taking $\varepsilon_0=\varepsilon_x^{-1}$. $\Box$

\subsection{Solutions of relative binomial Thue equations}

It follows from Theorem \ref{thbb} that in order to determine 
generators of relative power integral bases of $K$ over $L$ we have to
solve the relative binomial Thue equation (\ref{rrr}).

We have considered $-d=-3,-7,-11,-19,-43,-67,-163$. For all these 
values of $d$ we used the algorithm \cite{girelthue} to 
calculate the solutions with 
$\max(\overline{|X_0|},\overline{|Y_0|})<10^{250}$
of equation (\ref{rrr}) over $L=\Q(i\sqrt{d})$,
for all $m$ with
$1<m\leq 5000$, $m\equiv 2,3 \;(\bmod \ 4)$ and $(d,m)=1$.
In the following statements we list the solutions.

\begin{theorem}
Let $d$ be one of $-d=-3,-7,-11,-19,-43,-67,-163$, 
let $L=\Q(i\sqrt{d})$, $\omega=(1+i\sqrt{d})/2$.
For $1<m\leq 5000$, $m \equiv 2,3 (\bmod \ 4)$, $(d,m)=1$,
all solutions with
$\max(\overline{|X_0|},\overline{|Y_0|})<10^{250}$ of
equation (\ref{rrr}) in $(X_0,Y_0)\in \Z_L^2$
are \\
\\
{\bf 1.)} $(\pm 1,0)$ for all $d$ and $m$, \\
\\
{\bf 2.)} up to sign the following pairs for the listed $m$
for all $d$:\\
\noindent
$
m=2:\;\;\;\;  (1,1)\\
m=5 :\;\;\;\;  (3,2)\\
m=17 :\;\;\;\;  (2,1)\\
m=39 :\;\;\;\;  (5,2)\\
m=82 :\;\;\;\;  (3,1)\\
m=150:\;\;\;\;  (7,2)\\
m=257 :\;\;\;\;  (4,1)\\
m=410:\;\;\;\;  (9,2)\\
m=626 :\;\;\;\; (5,1)\\
m=915:\;\;\;\;  (11,2)\\
m=1297 :\;\;\;\; (6,1)\\
m=1785:\;\;\;\;  (13,2)\\
m=2402 :\;\;\;\;  (7,1)\\
m=3164:\;\;\;\;  (15,2)\\
m=4097:\;\;\;\;  (8,1)\\
$
\\
{\bf 3.)} up to sign the following pairs for the listed $m$ for special $d$: 
\\ \\
for $-d=-3$\\
\noindent
$
m=2:\;\;\;
                         (\omega,\omega),
                    (-\omega + 1,-\omega + 1),\\
m=5:\;\;\;
                       (3 \omega,2 \omega),
                  (-3 \omega + 3,-2 \omega + 2),\\
m=10:\;\;\;
                      (-2 \omega + 1,1),
                      (-\omega + 2,\omega),
                      (-\omega + 2,-\omega),
                    (\omega + 1,-\omega + 1),
                     (\omega + 1,\omega - 1),\\ \mbox{}\hspace{2cm}
                        (2 \omega - 1,1),\\
m=17:\;\;\;
                        (2 \omega,\omega),
                   (-2 \omega + 2,-\omega + 1),\\
m=39:\;\;\;
                       (5 \omega,2 \omega),
                  (-5 \omega + 5,-2 \omega + 2),\\
m=82:\;\;\;
                        (3 \omega,\omega),
                   (-3 \omega + 3,-\omega + 1),\\
m=145:\;\;\;
                       (-4 \omega + 2,1),
                     (-2 \omega + 4,\omega),
                    (-2 \omega + 4,-\omega),
                   (2 \omega + 2,-\omega + 1),\\ \mbox{}\hspace{2.2cm}
                   (2 \omega + 2,\omega - 1),
                       (4 \omega - 2,1),\\
m=150:\;\;\;
                      (7 \omega,2 \omega),
                 (-7 \omega + 7,-2 \omega + 2),\\
m=257:\;\;\;
                       (4 \omega,\omega),
                  (-4 \omega + 4,-\omega + 1),\\
m=410:\;\;\;
                      (9 \omega,2 \omega),
                 (-9 \omega + 9,-2 \omega + 2),\\
m=455:\;\;\;
                       (8,-2 \omega + 1),
                       (8,2 \omega - 1),
                     (8 \omega,-\omega + 2),
                   (-8 \omega + 8,\omega + 1),\\
m=580:\;\;\;
                      (17,-4 \omega + 2),
                       (17,4 \omega - 2),
                   (17 \omega,-2 \omega + 4),
                 (-17 \omega + 17,2 \omega + 2),\\
m=626:\;\;\;
                       (5 \omega,\omega),
                  (-5 \omega + 5,-\omega + 1),\\
m=730:\;\;\;
                       (-6 \omega + 3,1),
                     (-3 \omega + 6,\omega),
                    (-3 \omega + 6,-\omega),
                   (3 \omega + 3,-\omega + 1),\\ \mbox{}\hspace{2.2cm}
                   (3 \omega + 3,\omega - 1),
                       (6 \omega - 3,1),\\
m=905:\;\;\;
                      (19,-4 \omega + 2),
                       (19,4 \omega - 2),
                   (19 \omega,-2 \omega + 4),
                 (-19 \omega + 19,2 \omega + 2),\\
m=915:\;\;\;
                      (11 \omega,2 \omega),
                (-11 \omega + 11,-2 \omega + 2),\\
m=1111:\;\;\;
                      (10,-2 \omega + 1),
                      (10,2 \omega - 1),
                    (10 \omega,-\omega + 2),
                 (-10 \omega + 10,\omega + 1),\\
m=1297:\;\;\;
                       (6 \omega,\omega),
                  (-6 \omega + 6,-\omega + 1),\\
m=1785:\;\;\;
                     (13 \omega,2 \omega),
                (-13 \omega + 13,-2 \omega + 2),\\
m=2305:\;\;\;
                      (-8 \omega + 4,1),
                    (-4 \omega + 8,\omega),
                    (-4 \omega + 8,-\omega),
                  (4 \omega + 4,-\omega + 1),\\ \mbox{}\hspace{2.2cm}
                   (4 \omega + 4,\omega - 1),
                       (8 \omega - 4,1),\\
m=2402:\;\;\;
                       (7 \omega,\omega),
                  (-7 \omega + 7,-\omega + 1),\\
m=3164:\;\;\;
                     (15 \omega,2 \omega),
                (-15 \omega + 15,-2 \omega + 2),\\
m=4097:\;\;\;
                       (8 \omega,\omega),
                  (-8 \omega + 8,-\omega + 1),\\
$
\\
for $-d=-7$\\
\noindent
$
m=3:\;\;\;\;
                        (-2 \omega + 1,2),
                        (2 \omega - 1,2),\\
m=50:\;\;\;\;
                      (-2 \omega + 1,1),
                        (2 \omega - 1,1),\\
m=248:\;\;\;\;
                       (-6 \omega + 3,2),
                       (6 \omega - 3,2),\\
m=785:\;\;\;\;
                       (-4 \omega + 2,1),
                      (4 \omega - 2,1),\\
m=1914:\;\;\;\;
                      (-10 \omega + 5,2),
                      (10 \omega - 5,2),\\
m=3970:\;\;\;\;
                      (-6 \omega + 3,1),
                       (6 \omega - 3,1),\\
$
\\
for $-d=-11$\\
\noindent
$
m=122:\;\;\;\; (-2 \omega + 1,1), (2 \omega - 1,1)\\
m=1937 :\;\;\;\;(-4 \omega + 2,1), (4 \omega - 2,1)\\
$
\\
for $-d=-19$\\
\noindent
$
m=362 :\;\;\;\;(-2 \omega + 1,1), (2 \omega - 1,1)\\
$
\\
for $-d=-43$\\
\noindent
$
m=1850  :\;\;\;\;(-2 \omega + 1,1), (2 \omega - 1,1)\\
$
\\
for $-d=-67$\\
\noindent
$
m=4490 :\;\;\;\; (-2 \omega + 1,1) (2 \omega - 1,1)\\
$
\\
for $-d=-163$\\
\noindent
$
m=328 :\;\;\;\;(-2 \omega + 1,3), (2 \omega - 1,3).\\
$
\label{lista}
\end{theorem}

Using the solutions of the relative binomial Thue equations 
we construct the generators of relative power integral bases 
of $K$ over $L$.

\begin{theorem}
Let $d$ be one of $-d=-3,-7,-11,-19,-43,-67,-163$, 
let $L=\Q(i\sqrt{d})$, $\omega=(1+i\sqrt{d})/2$.
For $1<m\leq 5000$, $m \equiv 2,3 (\bmod \ 4)$, $(d,m)=1$,
all generators $\alpha=H+X\xi+Y\xi^2+Z\xi^3$ 
of relative power integral bases of $K$ over $L$ 
with $\max(\overline{|X|},\overline{|Y|},\overline{|Z|})<10^{500}$
are given by
\[
H\in \Z_L \; {\rm arbitrary},\;\;
X=X_0^2\varepsilon_0,\;\;
Y=\pm X_0Y_0\varepsilon_0,\;\;
Z=Y_0^2\varepsilon_0
\]
where $\varepsilon_0$ is a unit in $L$ and $(X_0,Y_0)$ is listed in 
Theorem \ref{lista}.
\end{theorem}

\noindent
{\bf Proof.} This is a consequence of Theorem \ref{thbb}. $\Box$

\subsection{Specialities of the actual calculations}

The algorithm for calculating small solutions of relative Thue equations
given in \cite{girelthue} consists of two parts. First the bound on the
variables (in our case $10^{250}$) is reduced by using LLL reduction algorithm.
In these examples the reduced bound was mostly between 10 and 200.
In the second step the tiny solutions under the reduced bound 
(say 200 in our case) are enumerated. In our present calculation
for relative binomial Thue equations
we reorganized these two parts to make the procedure more efficient.

We proceeded for the given values of $d$ separately.

The first part, the reduction was executed for all single $m$
with the fixed $d$. 
(This yields about 3000 values of $m$ allowed by the assumptions.)
These (about 3000) reduction procedures were executed (all together) within about 
3.2--3.6 hours of CPU time for each $d$. 

The second part of the procedure described in \cite{girelthue}, 
the enumeration of tiny values of the variables
(with absolute values under the reduced bound, say 200 in our examples)
was performed a different way. 
We proceeded for all values of $d$ separately.
It turned out to be very CPU time
consuming to test the tiny values of $x_1,y_1,x_2,y_2$ for all
$m$ if they satisfy the relative Thue equation. On the other hand
we have run the cycles for all $x_1,y_1,x_2,y_2$ under the reduced 
bound and determined those $m$ for which they yield a
solution. Using integer arithmetic and obvious refinements
in the enumeration process, to test all $x_1,y_1,x_2,y_2$
and finding the suitable values of $m$ took about 3.5 hours.

\subsection{CPU times}

The CPU times of the paper refer to an average laptop. The
programs were developed in Maple \cite{maple} and were executed under Linux.
The same programs were also tested 
on the high performance computer (supercompuer)
of the University of Debrecen
where we obtained about 20\% better results, calculated for a single
processor node.

\end{document}